\documentclass[10pt,english]{smfart}
\usepackage{amssymb}
\usepackage{amsmath}
\usepackage{amscd}
\DeclareFontEncoding{OT2}{}{} 
  \newcommand{\textcyr}[1]{%
    {\fontencoding{OT2}\fontfamily{wncyr}\fontseries{m}\fontshape{n}%
     \selectfont #1}}
\newcommand{\Sha}{{\mbox{\textcyr{Sh}}}}
\font\notsosmall=ptmr8r scaled 450
\font\small=ptmr8r scaled 550
\font\realsmall=ptmr8r scaled 400
\def\P{{\bf P}}
\def\Z{{\bf Z}}
\def\Cl{{\mathcal Cl}}
\def\C{{\bf C}}
\def\R{{\bf R}}
\def\Q{{\bf Q}}
\def\WP{{\mathcal P}}
\def\SS{{\mathcal S}}
\def\H{{\mathcal H}}
\def\HH{{\bf H}}
\def\BAR{\overline}
\def\Re{{\rm Re}}
\def\Im{{\rm Im}}
\usepackage{amssymb}
\usepackage{amsmath}
\usepackage{amscd}
\theoremstyle{plain}
\newtheorem{prop}{Proposition}[section]
\newtheorem{proposition}[prop]{Proposition}
\newtheorem{theorem}[prop]{Theorem}
\newtheorem{conjecture}[prop]{Conjecture}
\newtheorem{algorithm}[prop]{Algorithm}
\newtheorem{subalgorithm}[prop]{Subalgorithm}
\newtheorem{lemma}[prop]{Lemma}
\theoremstyle{definition}
\newtheorem{definition}[prop]{Definition}
\def\kill#1{}%
\begin{document}
\title{Some remarks on Heegner point computations}
\alttitle{Une partie remarque sur des calculs de point de Heegner}
\author{Mark Watkins}
\address
{School of Mathematics, University Walk,
University of Bristol, Bristol BS8 1TW, ENGLAND}
\email{watkins@maths.usyd.edu.au}
\begin{abstract}
We give an overview of the theory of Heegner points for elliptic curves,
and then describe various new ideas that can be used in the computation
of rational points on rank 1 elliptic curves.
In particular, we discuss the idea of Cremona (following Silverman)
regarding recovery a rational point via knowledge of its height, the idea of
Delaunay regarding the use of Atkin-Lehner involutions in the selection of
auxiliary parameters, and the idea of Elkies regarding descent and lattice
reduction that can result in a large reduction in the needed amount of
real-number precision used in the computation.
\end{abstract}
\kill{
\begin{altabstract}
Nous donnons une vue d'ensemble de la th\'eorie de points de Heegner
pour les courbes elliptiques, et puis d\'ecrivons les diverses
nouvelles id\'ees qui peuvent \^etre empoly\`ees
dans le calcul des points raisonnables sur les
courbes elliptiques du grade 1.
En particulier, nous discutons l'id\'ee de Cremona
(apr\`es Silverman) concernant le r\'etablissement par point
raisonnable par l'interm\'ediaire de la connaissance de sa taille,
l'id\'ee de Delaunay concernant l'utilisation des involutions d'Atkin-Lehner
dans le choix des param\`etres auxiliaires, et l'id\'ee 'Elkies
concernant la descente et la r\`eduction de trellis qui peuvent avoir comme
cons\`equence une grande r\'eduction de la quantit\'e n\'ecessaire
de pr\'ecision de vrai-nombre utilis\'ee dans le calcul.
\end{altabstract}}
\altkeywords{elliptic curves, Heegner points, descent, lattice reduction}
\maketitle

\section{Introduction}
We make some remarks concerning Heegner point computations.
One of our goals shall be to give an algorithm (perhaps conditional
on various conjectures) to find a non-torsion rational point
on a given rank 1 elliptic curve.
Much of this is taken from a section in Henri Cohen's
latest book~\cite{Cohen}, and owes a great debt to Christophe Delaunay.
The ideas in the section about lattice reduction
are largely due to Noam Elkies.
We do not delve deeply into the theory of Heegner points, but simply
give references where appropriate; the recent MSRI book
``Heegner points and Rankin $L$-series'' \cite{MSRI} contains many
good articles which consider Heegner points and their generalisations
from the standpoint of representation theory.

The author thanks the Institut Henri Poincar\'e for its hospitality
and the Centre National de la Recherche Scientifique for financial support.
The algorithms described here have been implemented in the Magma computer
algebra system \cite{magma};
the author thanks the Magma computer algebra group at
the University of Sydney for their hospitality and financial support.
The author was also partially funded by an NSF VIGRE Postdoctoral
Fellowship during part of this work.

\section{Definitions and Outline of Theory}
Let $\tau$ be a quadratic surd in the upper-half-plane~$\HH$.
Let $f_\tau=(A,B,C)$ be the associated integral primitive
positive-definite binary quadratic form,
so that \hbox{$A\tau^2+B\tau+C=0$} with $A>0$
and~$\gcd(A,B,C)=1$. The discriminant $\Delta(\tau)$
is $\Delta(f_\tau)=B^2-4AC$, which is negative.
For simplicity we take $\Delta(\tau)$ to be fundamental,
though much of the theory can be made to work when it is not.
However, consideration of positive discriminant does not follow
in the same manner.
\begin{definition}
A Heegner point of level~$N$ and discriminant~$D$ is a quadratic surd
in the upper-half-plane with $\Delta(\tau)=D=\Delta(N\tau)$.
We let $\H_N^D$ be the set of Heegner points of level~$N$ and
discriminant~$D$.
\end{definition}
\begin{proposition}
Let $\tau\in\HH$ be a quadratic surd with discriminant~$D$
and~$f_\tau=(A,B,C)$. Then $\tau\in\H_N^D$ iff $N|A$ and $\gcd(A/N,B,CN)=1$.
\end{proposition}
\begin{proof}
Note that $\tau={-B+\sqrt D\over 2A}$ and $N\tau={-NB+N\sqrt D\over 2A}$.
For $\Delta(\tau)=\Delta(N\tau)$ we need $N\tau={-B'+\sqrt D\over 2A'}$
and by equating imaginary and real parts we get $A=NA'$ and~$B=B'$,
so that $N|A$. Also note that $(A/N)(N\tau)^2+B(N\tau)+(CN)=0$,
from which we get the rest of the lemma.
\end{proof}
Note that $\H_N^D$ will be empty unless $D=B^2-4N(A/N)C$
is a square modulo~$4N$.
\begin{lemma}
The set $\H_N^D$ is closed under $\Gamma_0(N)$-action.
\end{lemma}
\begin{proof}
If $\gamma\in SL_2(\Z)\supseteq \Gamma_0(N)$ then
$\Delta\bigl(\gamma(\tau)\bigr)=\Delta(\tau)$ since
the discriminant is fixed. A computation shows that
$\gamma\in\Gamma_0(N)$ and $\tau\in\H_N^D$ imply~$\gamma(\tau)\in\H_N^D$.
\end{proof}
\begin{lemma}
The set $\H_N^D$ is closed under the $W_N$-action that
sends $\tau\rightarrow -1/N\tau$.
\end{lemma}
\begin{proof}
Follows from the above proposition since $f_{(-1/N\tau)}=(CN,-B,A/N)$.
\end{proof}
\begin{definition}
We let $\SS(D,N)$ be the set of square roots mod~$2N$ of $D$ mod~$4N$.
\end{definition}
\begin{theorem}
The sets $\H_N^D/\Gamma_0(N)$ and
$\SS(D,N)\times \Cl\bigl(\Q(\sqrt D)\bigr)$ are in bijection.
\end{theorem}
\begin{proof}
This can be shown by chasing definitions.
Essentially, $[\tau]\in\H_N^D/\Gamma_0(N)$ gets mapped to
$(B$ mod $2N)\times [\Z+\tau\Z]$ where $f_\tau=(A,B,C)$,
and in the other direction,
when we are given $\beta\times l\in \SS(D,N)\times \Cl\bigl(\Q(\sqrt D)\bigr)$
we take $(A,B,C)\in l$ with $N|A$ and $B\equiv\beta$ (mod~$2N$),
and then~$\tau={-B+\sqrt D\over 2A}$.
\end{proof}
From now on we let $E$ be a global minimal model of a rational elliptic curve
of conductor~$N$, and take $D$ to be a negative fundamental discriminant such
that $D$ is a square mod $4N$.
We let $\HH^\star$ be the union of $\HH$ with the rationals and~$i\infty$.
We let $\WP(z)$ be the function that sends $z\in\C/\Lambda$ to the point
$\bigl(\wp(z),\wp^\prime(z)\bigr)$ on~$E$.
\begin{theorem}
There is a surjective map \hbox{$\hat\phi:X_0(N)\rightarrow E$}
(the modular parametrisation)
where $X_0(N)=\HH^\star/\Gamma_0(N)$ and $E$ can be viewed
as $\C/\Lambda$ for some lattice~$\Lambda$.
This map can be defined over the rationals.
\end{theorem}
\begin{proof}
This is due to Wiles and others \cite{wiles,TW,diamond,CDT,BCDT}.
We let $\phi$ be the associated map from $\HH^\star/\Gamma_0(N)$
to~$\C/\Lambda$.
Explicitly, we have that $\tau\in\HH^\star$ gets mapped to the complex point
$\phi(\tau)=2\pi i\int_{i\infty}^\tau \psi_E=\sum_n (a_n/n)e^{2\pi i n\tau}$,
where $\psi_E$ is the modular form of weight 2 and level~$N$ associated to~$E$.
The lattice~$\Lambda$ is generated by the
real and imaginary periods,\footnote{
 Our convention is that the imaginary period is purely imaginary
 when the discriminant of~$E$ is positive, and in the negative discriminant
 case the real part of the imaginary period is~$\Omega_{\rm re}/2$.
 The fundamental volume~$\Omega_{\rm vol}$
 is the area of the period parallelogram.}
which we denote by $\Omega_{\rm re}$ and~$\Omega_{\rm im}$.
We assume that the Manin constant is~1,
which is conjectured always to be the case for curves of positive rank
(see~\cite{Stevens} and~\cite{SteinWatkins}).
\end{proof}
\begin{theorem}\label{thm:cm}
Let $\tau=\beta\times l\in\H_N^D$.
Then $\WP\bigl(\phi(\tau)\bigr)$ has its coordinates
in the Hilbert class field of~$\Q(\sqrt D)$.
Also we have
\begin{enumerate}
\item $\BAR{\phi(\beta\times l)}=\phi\bigl(-\beta\times l^{-1}\bigr)$,
in~$\C/\Lambda$.
\item $\phi\bigl(W_N(\beta\times l)\bigr)=
\phi\bigl(-\beta\times ln^{-1}\bigr)$ in $\C/\Lambda$ where
$n=\bigl[ N\Z+{\beta+\sqrt D\over 2}\Z\bigr]$,
\item
$\WP\bigl(\phi(\beta\times l)\bigr)^{{\rm Artin}(m)}
=\WP\bigl(\phi\bigl(\beta\times lm^{-1}\bigr)\bigr)$
for all $m\in\Cl\bigl(\Q(\sqrt D)\bigr)$.
\end{enumerate}
\end{theorem}
\begin{proof}
This is the theorem of complex multiplication
of Shimura~\cite{ShimuraTaniyama,Shimura}.
We outline the proof of the first statement,
for which we work via the modular $j$-function.
We have that $j(\tau)$ is in the Hilbert class field $H$
(see \cite[II, 4.3]{silverman2}) and similarly with~$j(N\tau)$.
Thus we get that $X_0(N)$ over $H$ contains the moduli point
corresponding to the isogeny between curves with these $j$-invariants.
Since the modular parametrisation map $\hat\phi$ can be defined over the
rationals, the image of the moduli point under $\hat\phi$ has its
coordinates in the Hilbert class field.
\end{proof}
Note that $\WP\bigl(\overline{\phi(\tau)}\bigr)=
\overline{\WP\bigl(\phi(\tau)\bigr)}$, so that there
is no danger of confusing complex conjugation in~$\C/\Lambda$
with complex conjugation of the coordinates of the point on~$E$.
Using the third fact of Theorem~\ref{thm:cm}, we can take the trace
of~$\WP\bigl(\phi(\tau)\bigr)$
and get a point that has coordinates in~$\Q(\sqrt D)$.
Indeed, writing $H$ for the Hilbert class field and $K$ for
the imaginary quadratic field $\Q(\sqrt D)$ we get that
\begin{align*}
P={\rm Trace}_{H/K}\bigl(\WP\bigl(\phi(\tau)\bigr)\bigr)
&=\sum_{\sigma\in{\rm Gal}(H/K)} \WP\bigl(\phi(\tau)\bigr)^\sigma=
\sum_{m\in\Cl(K)} \WP\bigl(\phi(\beta\times l)\bigr)^{{\rm Artin}(m)}\\
&\qquad=\sum_{m\in\Cl(K)} \WP\bigl(\phi\bigl(\beta\times lm^{-1}\bigr)\bigr)=
\sum_{m\in\Cl(K)} \WP\bigl(\phi(\beta\times m)\bigr)
\end{align*}
has coordinates in~$\Q(\sqrt D)$. When $E$ has odd functional equation,
we can use the first two facts of Theorem~\ref{thm:cm} to show that~$P=\BAR P$,
so that $P$ has coordinates in~$\Q$. In this case we have
$\psi_E=\psi_E \circ W_N$ which implies~$\phi=\phi\circ W_N$,
so that in $\C/\Lambda$ we have
$$\BAR{\phi(\beta\times m)}=\BAR{\phi\bigl(W_N(\beta\times m)\bigr)}
=\BAR{\phi\bigl(-\beta\times mn^{-1}\bigr)}=
\phi\bigl(\beta\times m^{-1}n\bigr),$$
which gives us that
$$\BAR P=\sum_{m\in\Cl(K)} \WP\bigl(\BAR{\phi(\beta\times m)}\bigr)
=\sum_{m\in\Cl(K)} \WP\bigl(\phi\bigl(\beta\times m^{-1}n\bigr)\bigr)
=\sum_{m\in\Cl(K)} \WP\bigl(\phi(\beta\times m)\bigr)=P.$$
We can rewrite some of this by introducing some new notation.
\begin{definition}
We write $\H_N^D(\beta)$ for subset of $\tau\in\H_N^D$ such that the
associated form $f_\tau=(A,B,C)$ has $B\equiv\beta$ (mod~$2N$).
We write $\hat\H_N^D=\H_N^D/\Gamma_0(N)$, and noting that $\Gamma_0(N)$ acts
on $\H_N^D(\beta)$, we write $\hat\H_N^D(\beta)=\H_N^D(\beta)/\Gamma_0(N)$.
\end{definition}
Since $\hat\H_N^D(\beta)$ is in 1-1 correspondence
with $\Cl\bigl(\Q(\sqrt D)\bigr)$, we get that
$$P=\sum_{m\in\Cl(Q(\sqrt D))} \WP\bigl(\phi(\beta\times m)\bigr)=
\sum_{\tau\in\hat\H_N^D(\beta)} \WP\bigl(\phi(\tau)\bigr).$$

\section{The Gross-Zagier theorem and an algorithm}
We now have a plan of how to find a non-torsion point
on a curve of analytic rank~1.
We select an auxiliary negative fundamental discriminant~$D$ such that
$D$ is a square modulo~$4N$, choose~$\beta\in \SS(D,N)$,
find $\tau$-representatives for~$\hat\H_N^D(\beta)$,
compute $\phi(\tau)$ for each, sum these in $\C/\Lambda$,
map the resulting point to $E$ via the Weierstrass parametrization,
and try to recognize the result as a rational point. One problem
is that we might get a torsion point. Another problem is that we won't
necessarily get a generator, and thus the point might have inflated height,
which would increase our requirements on real-number precision.
The Gross-Zagier Theorem tells us what height to expect, and combined with
the Birch--Swinnerton-Dyer Conjecture, we get a prediction of what height
a generator should have. Our heights will be the ``larger'' ones, and are
thus twice those chosen by some authors.
\begin{theorem}
Suppose $D<-3$ is a fundamental discriminant
with $D$ a square modulo~$4N$ and $\gcd(D,2N)=1$. Then
$$h(P)={\sqrt{|D|}\over 4\Omega_{\rm vol}}L'(E,1)L(E_D,1)\times
2^{\omega(\gcd(D,N))}\biggl({w(D)\over 2}\biggr)^2.$$
\end{theorem}
\begin{proof}
This is due to Gross and Zagier~\cite{GrossZagier}.
Here $E_D$ is the quadratic twist of~$E$ by~$D$,
while $w(D)$ is the number of units in $\Q(\sqrt D)$
and $\omega(n)$ is the number of distinct prime factors of~$n$.
\end{proof}
Calculations of Gross and Hayashi~\cite{Hayashi} indicate that
this height formula is likely to be true for all negative fundamental
discriminants~$D$ that are square mod~$4N$.

We now write $P=lG+T$ where $G$ is a generator\footnote
{Note that we will actually get $\sqrt{\#\Sha}$ times a generator
 from our method, since we cannot disassociate $\Sha$ from the
 regulator in the Birch--Swinnerton-Dyer formula.}
and $T$ is a torsion point,
so that $h(P)=l^2h(G)$. Then we replace $L'(E,1)$ through use of the
Birch--Swinnerton-Dyer conjecture~\cite{BSD} to get the following:\footnote{
 We use the convention that the Tamagawa number at infinity is equal to the
 number of connected components of~$E$ over~$\R$ --- thus it is 1 for curves
 with negative discriminant and 2 for curves with positive discriminant.}
\begin{conjecture}
With notations as above we have that
$$l^2={\Omega_{\rm re}\over 4\Omega_{\rm vol}}
\biggl({\prod_{p|N\infty} c_p\cdot\#\Sha}\biggr)
{\sqrt{|D|}\over \#E(\Q)_{\rm tors}^2}L(E_D,1)\cdot
\biggl({w(D)\over 2}\biggr)^2 2^{\omega(\gcd(D,N))}.$$
\end{conjecture}
In particular, we note that we should use a quadratic twist~$E_D$
that has rank zero, so that $L(E_D,1)$ does not vanish.
The existence of such a twist is proven in~\cite{BumpFriedbergHoffstein}.
Thus we have the following algorithm, which we shall work on improving.
\begin{algorithm}
Given a rational elliptic curve~$E$ of conductor~$N$ of analytic rank~1,
find a non-torsion rational point.
\begin{enumerate}
\item Compute $L'(E,1)$ and find a fundamental discriminant $D<0$ with $D$
a square modulo~$4N$ and $L(E_D,1)\neq 0$, so that the index~$l$ is nonzero.
\item Choose $\beta\in \SS(D,N)$ and compute (to sufficient precision)
the complex number
$$z=\sum_{\tau\in\hat\H_N^D(\beta)}\phi(\tau)=
\sum_{\tau\in\hat\H_N^D}\sum_{n=1}^\infty {a_n\over n}e^{2\pi i n\tau}.$$
\item Let $m$ be the gcd of~$l$ and the exponent of the torsion group of~$E$.
If the discriminant of $E$ is positive, check if $\WP(\dot z)$
is close to a rational point on~$E$ for $u=1,\dots,lm$ for both
$$\dot z=(m\Re(z)+u\Omega_{\rm re})/ml\quad{\rm and}\quad
\dot z=(m\Re(z)+u\Omega_{\rm re})/ml+\Omega_{\rm im}/2.$$
\vskip0pt\noindent
If the discriminant of $E$ is negative,
let $o=\Im(z)/\Im(\Omega_{\rm im})$ and check $\WP(\dot z)$ for
$\dot z=(m\Re(z)+u\Omega_{\rm re})/ml+o\Omega_{\rm re}/2$
over the same \hbox{$u$-range}.
\end{enumerate}
\end{algorithm}
One can compute the index~$l$ in parallel with the~$\phi(\tau)$,
since both involve computing the $a_n$ of the elliptic curve~$E$.
However, this can cause problems if the index turns out to be zero
(that is, if $E_D$ has positive rank).

\subsection{Step 2 of the Algorithm}
We now discuss how to do the second step efficiently.
First note that we can sometimes pair $\phi(\tau)$ with its complex conjugate;
recalling that $\phi=\phi\circ W_N$,
by Theorem~\ref{thm:cm} in $\C/\Lambda$ we have
\begin{equation}\label{eq:chap3}
\overline{\phi(\beta\times l)}=\phi(-\beta\times l^{-1})=
\phi\bigl(W_N(-\beta\times l^{-1})\bigr)=
\phi\bigl(\beta\times(ln)^{-1}\bigr).
\end{equation}
For $f=(A,B,C)\in\hat\H_N^D$ we write~$\BAR f=(A/N,-B,CN)$,
so that when~$g\sim\BAR f$ in the class group we have $\phi(f)=\BAR{\phi(g)}$
and thus $\phi(f)+\phi(g)=2\,\Re\,\phi(f)$ in~$\C/\Lambda$.
We refer to this as {\bf pairing} the forms.

For $\sum_n (a_n/n)e^{2\pi i n\tau}$ to converge rapidly, we wish for
the imaginary parts of our representative~$\tau$'s to be large.
It turns out the best we can do is essentially have the smallest
imaginary part be about $1/N$ in size. We can achieve this via
a trick of Delaunay, which introduces more Atkin-Lehner involutions.
\begin{definition}
Let $Q|N$ with $\gcd(Q,N/Q)=1$, and let $u,v\in\Z$ be such that $uQ^2-vN=Q$.
The Atkin-Lehner involution $W_Q$ sends $\tau$ to ${uQ\tau+v\over N\tau+Q}$.
\end{definition}
This defines $W_Q$ up to transformations by elements in~$\Gamma_0(N)$.
One can check that $W_Q(W_Q(\tau))$ is in the $\Gamma_0(N)$-orbit of~$\tau$,
and that the $W_Q$ form an elementary abelian 2-group~$W$
of order~$2^{\omega(N)}$.
The important fact about the $W_Q$'s that we shall use is that
\hbox{$\psi_E=\pm\psi_E\circ W_Q$}, so that
$\phi(\tau)=\pm\phi(W_Q(\tau))+\phi\bigl(W_Q(i\infty)\bigr)$.
The sign can be computed as
$\epsilon_Q=\prod_{p|Q} \epsilon_p$ where $\epsilon_p$
is the local root number of~$E$ at~$p$. Delaunay's idea is to maximise
the imaginary part of~$\tau$ over $\Gamma_0(N)$ and $W$ rather than
just $\Gamma_0(N)$; the difficulty is that the action of $W_Q$
need not preserve~$\beta$. However, we still have that
$$P=\sum_{\tau\in\hat\H_N^D(\beta)} \phi(\tau)=
\sum_{\tau\in\hat\H_N^D(\beta)} \epsilon_Q\phi\bigl(W_Q(\tau)\bigr)
+({\rm torsion\>point}).$$
For the analogue of the second part of Theorem~\ref{thm:cm},
we need to consider what happens with to~$\beta$.
We define $\beta_Q$ as follows.
We make $\beta_Q$ and $\beta$ have opposite signs mod~$p^k$ for prime
powers~$p^k$ with $p^k\|Q$ and~$\gcd(p,D)=1$, and else $\beta=\beta_Q$.
In particular, we have that~$Q=\gcd(\beta-\beta_Q,N)$ when~$N$ is odd.
The desired analogue is now that
$$\phi\bigl(W_Q(\beta\times l)\bigr)=
\epsilon_Q\phi(\beta_Q\times lq^{-1})+
\phi\bigl(W_Q(i\infty)\bigr)\quad{\rm with}
\quad q=\biggl[Q\Z+{-\beta_Q+\sqrt D\over 2}\Z\biggr].$$

The primes $p$ which divide $D$ are different since there is
only one square root of~$D$ mod~$p$; thus $\beta$ is preserved upon
application of~$W_Q$ for $Q$ that are products of such primes.
For such~$Q$, we can note the following with respect to complex conjugation.
Suppose we have that $m\sim (ln)^{-1}$ so that by \eqref{eq:chap3} we have
$\BAR{\phi(\beta\times l)}=\phi(\beta\times m)$.
Then, using the fact that $q^{-1}=q$ in this case,
in $\C/\Lambda$ we have, up to torsion, that
\begin{align*}
\BAR{\phi\bigl(W_Q(\beta\times l)\bigr)}&\circeq
\epsilon_Q\BAR{\phi(\beta_Q\times lq^{-1})}=
\epsilon_Q\phi\bigl(\beta_Q\times (lq^{-1}n)^{-1}\bigr)=
\epsilon_Q\phi\bigl(\beta_Q\times (lqn)^{-1}\bigr)\circeq\\
&\circeq\phi\bigl(W_Q(\beta\times (ln)^{-1})\bigr)=
\phi\bigl(W_Q(\beta\times m)\bigr).
\end{align*}
So we see that $(\beta\times l)$ can be paired iff $W_Q(\beta\times l)$
can be paired.

We now give the algorithm for finding good~\hbox{$\tau$-representatives.}
The idea to run over all forms $(aN,b,c)$ of discriminant~$D$ with
$a$ small, mapping these via the appropriate Atkin-Lehner involution(s)
to forms with fixed square root~$\beta$, and doing this until the images
cover the class group. Of course, the conjugation action is also considered.
\begin{subalgorithm}
Given $D,N$, find good $\tau$-representatives.
\begin{enumerate}
\item Choose $\beta\in \SS(D,N)$. Set $U=\emptyset$ and $R=\emptyset$.
\item While $\#R\neq\#\Cl(\bigl(Q\sqrt D)\bigr)$ do:
\item \quad Loop over $a$ from 1 to infinity and~$b\in \SS(D,N)$
[lift $b$ from $\Z/2N$ to~$\Z$]:
\item \quad\quad Loop over all solutions~$s$
of $Ns^2+bs+(b^2-D)/4N\equiv 0$ modulo~$a$:
\item \quad\quad\quad Let $f=\bigl(aN,b+2Ns,((b+2Ns)^2-D)/4aN\bigr)$.
\item \quad\quad\quad Loop over all positive divisors~$d$ of~$\gcd(D,N)$
[which is squarefree]:
\item \quad\quad\quad\quad Let $g=W_Q(f)/Q$ where~$Q$
is $d$ times the product of the $p^k\|N$
\break\hphantom0\quad\quad\quad\quad
with \hbox{$b\not\equiv\beta$ mod~$p^k$,}
so that $g\in \H_N^D(\beta)$.
\item \quad\quad\quad\quad
If the reductions of $g$ and $\BAR g$ are both not in~$R$ then
append them\break\hphantom0\quad\quad\quad\quad
to~$R$, and append $f$ to $U$ with weight $\epsilon_Q$
when $g\sim\BAR g$ and with\break\hphantom0\quad\quad\quad\quad
weight $2\epsilon_Q$ when~$g\not\sim\BAR g$.
\end{enumerate}
\end{subalgorithm}
With this subalgorithm,
we get that~$z=\sum_{f\in U} {\rm weight}(f)\phi(\tau_f)$
in Step 2 of the main algorithm.
We expect the maximal~$a$ to be of size $\#\Cl(\Q(\sqrt D))/2\#W$.
This subalgorithm makes
``parameter selection'' fast compared to the computation of the~$\phi(\tau)$.
\subsection{Step 3 of the Algorithm}
We now turn to the last step of our main algorithm,
reconstructing a rational point
on an elliptic curve from a real approximation. The most na\"{\i}ve method
for this is simply to try to recognise the $x$-coordinate as a rational
number. If our height calculation tells us to expect a point whose
\hbox{$x$-coordinate} has a numerator and denominator of about $H$ digits,
the use of continued fractions will recognise it if we do all computations
to about twice the precision, or $2H$~digits. We can note that by using
a \hbox{degree-$n$} map to $\P^{n-1}$ and \hbox{$n$-dimensional} lattice
reduction, this can be reduced to $nH/(n-1)$ digits
for every~$n\ge 3$ --- we will discuss a similar idea later when
we consider combining descent with our Heegner point computations.
But in this case we can do better;
we are able to recognise our rational point with only $H$ digits of
precision due to a trick of Cremona, coming from an idea in a paper
of Silverman~\cite{Silverman}.
The idea is that we know the canonical height of our desired point,
and this height decomposes into local heights; we have
$$h(P)=h_{\infty}(P)+\sum_{p|N} h_p(P)+\log{\rm denominator}(x(P)).$$
The height at infinity $h_{\infty}(P)$ can be approximated from a
real-number approximation to~$P$, and there are finitely many possibilities
for each local height $h_p(P)$ depending the reduction type of~$E$ at~$p$.
We compute the various local heights to $H$ digits of precision, and then
can determine the denominator of $x(P)$ from this, our task being eased
from the fact that it is square. Then from our real-approximation\footnote{
 Elkies tells us that, given the height to precision~$H$,
 the techniques of \cite{elkies-ants4} (see Theorem~4 in particular)
 can reduce the needed precision of the real Heegner approximation
 to $o(H)$ as~$H\rightarrow\infty$.
 The idea is that for a fixed~$C$ the equation
 \hbox{$h_{\infty}\bigl((x,y,z)\bigr)+2\log z=C$}
 defines a {\it transcendental} arc, and thus the use of a sufficiently high
 degree Veronese embedding will reduce the needed precision substantially.
 This method in its entirety might not be that practical, though the
 use of height information in conjunction with the geometry of the curve
 should allow a useful reduction in precision.}
of~$P$ we can recover the \hbox{$x$-coordinate}, and from this
we get~$P$. Note that we need to compute $L'(E,1)$ to a precision
of $H$ digits, but this takes only about $\sqrt N(\log H)^{O(1)}$ time.
In practise, there can be many choices for the sums of local heights,
and if additionally the index is large, then this step can be quite
time-consuming. This can be curtailed a bit by doing the calculations for
the square root of the denominator of the $x$-coordinate to only
about $H/2$~digits, and then not bothering with the elliptic exponential
step unless the result is sufficiently close to an integer.

\subsection{Example}
We now give a complete example. Other explicit descriptions of computations
with Heegner points appear in \cite{BirchStephens,Stephens,Liverance}.
We take the curve given by $[1,-1,0,-751055859,-7922219731979]$
for which the Heegner point has height $139.1747+$. We select $D=-932$,
for which the class number is~12 and the index $l$ is~4.
We have $N=11682$ and choose~$\beta=214$.
Our first form is $(11682,214,1)$ to which we apply~$W_1={\rm id}$.
The reduction of this is $(1,0,233)$, and it pairs with itself
under complex conjugation. Since we have $\gcd(D,N)=2$, we can use~$W_2$
without changing~$\beta$; we get the form
$(206717861394,70769770,6057)$ which reduces
to the self-paired form~$(2,2,117)$.
Our next form is $(11682,2338,117)$ to which we apply~$W_{11}$
to get $(122225810454,230158978,108351)$ which reduces to
$(11,6,22)$ and pairs with~$(11,-6,22)$. Applying~$W_{22}$ gives
a form which reduces to~$(11,-6,22)$, so we ignore it.
Next we have $(11682,2810,169)$ to which we apply~$W_9$,
getting a form that reduces to $(9,2,26)$ and pairs with~$(9,-2,26)$.
Applying~$W_{18}$ gives a form that reduces to $(13,-2,18)$
and pairs with~$(13,2,18)$.
Then we have $(11682,4934,521)$ to which we apply~$W_{99}$,
getting a form that reduces to $(3,-2,78)$ and pairs with~$(3,2,78)$.
And finally applying $W_{198}$ we get a form
that reduces to $(6,-2,39)$ and pairs with~$(6,2,39)$,
and so we have all of our~\hbox{$\tau$-representatives}.
We note that $W_{11}$, $W_9$, and $W_{18}$ switch the sign of the
modular form, and thus the obtained forms get a weighting of~$-2$.
The self-paired forms get a weighting of~$+1$, and the other two forms
get a weighting of~$+2$. For the non-self-paired forms we must remember
to take the real part of the computed~$\phi(\tau)$
when we double it.\footnote{
 The self-paired forms~$f$ have $\phi(f)=\BAR{\phi(f)}$ in $\C/\Lambda$
 but not necessarily in~$\C$ --- the imaginary part cannot be ignored
 when the discriminant of~$E$ is negative and $lm$ is odd.}
The pairing turns is rather simple in this example,
but need not be so perspicacious with respect to the class group.
Note that we use only four distinct forms for our computations.

We need about 60 digits of precision if we use the Cremona-Silverman
method to reconstruct the rational point, which means we must compute
about 20000 terms of the~\hbox{$L$-series}. The curve $E$ has negative
discriminant and no rational torsion points. We compute a real-approximation
to the Heegner point in $\C/\Lambda$ to be
$$z=0.00680702983101357730368201485198918786991251635619740952608094.$$
We have $o=\Im(z)/\Im(\Omega_{\rm im})=0$, and with $l=4$ and $u=2$ we get that
$$\dot z=0.00891152819280235244790996808333469812474933020620405901507952,$$
to which we apply the Cremona-Silverman method of recovery.
The curve~$E$ is annoying for this method,
in that we have many possibilities for~$h_p(P)$.
The height of the Heegner point is given by
$$h(P)=139.174739524758127811521877478222781093487974225206369462318,$$
and the height at infinity\footnote{
 Note that Silverman \cite{Silverman}
 uses a different normalisation of height,
 and his choice of the parameter~$z$ when he computes the height at infinity
 corresponds to $\dot z/\Omega_{\rm re}$ for us.
 Also, his method is only linearly convergent,
 while that given in \cite[\S 7.5.7]{cohen0} is quadratically convergent.}
is given by
$$h_\infty(P)=2.10306651755149369196435189022120441716979687181328497567075.$$
The reduction type at~2 is~$I_{25}$, at~3 it is $I_{13}^\star$,
at~11 it is~$I_1$, and at~59 it is~$I_3$.
Thus we have $13\times 3\times 1\times 2$ choices for the local heights. 
It turns out that we have\footnote{
 This follows {\it a posterori} since $P$ is nonsingular
 modulo these primes of multiplicative reduction.}
$h_p(P)={1\over 6}v_p(\Delta)\log p$ for $p=2,11,59$,
while~$h_3(P)=(13/6)\log 3$.
The denominator of the $x$-coordinate is
$12337088946900997614694947283^2$, and the numerator is
$$5908330434812036124963415912002702659341205917464938175508715.$$

\subsection{Variants}
Next we mention a variant which, for the congruent number curve,
has been investigated in depth by Elkies~\cite{Elkies}.
Here we fix a rank zero curve, say the curve $E: y^2=x^3-x$
of conductor~32, and try to find points on rank~1 quadratic twists~$E_D$
with $D<0$. It can be shown that $E_D$ will have odd functional equation for
\hbox{$|D|\equiv 5,6,7$ (mod~8).} There is not necessarily a Gross-Zagier
theorem in all these cases, and some involve mock Heegner points instead
of Heegner points. However, we still have the prediction that
$h(P)=\alpha_D L(E,1) L'(E_D,1)$ for some~\hbox{$\alpha_D>0$.}
Elkies computes a point $P$ in $\C/\Lambda$ via a method
similar to the above --- however, he generally\footnote{
 Elkies also computes a generator of height $239.6+$
 for $y^2=x^3-1063^2x$ in this manner.}
only attempts to determine if it is non-torsion,
and thus need not worry as much about precision.
There are about $\#\Cl\bigl(\Q(\sqrt D)\bigr)\approx\sqrt{|D|}$
conjugates of~$\tau$ for which $\phi(\tau)$ needs to be computed;
since we have an action of~$\Gamma_0(32)$, computing each $\phi(\tau)$
takes essentially constant time, so we get an algorithm that takes
about time $|D|^{1/2}$ to determine whether the computed point is non-torsion.
Note that we don't obtain $L'(E_D,1)$, which takes about $|D|$ time to
compute, but only whether it is nonzero.
MacLeod~\cite{MacLeod} investigated a similar family of quadratic twists,
those of a curve of conductor~128.
The relevant curves are $y^2=(x+p)(x^2+p^2)$ with $p\equiv 7$ (mod~8);
with $p=3167$ the height is $1022.64+$.
Some additional papers that deal with the theory and constructions in
this case are those of Birch and Monsky \cite{birch,birch2,monsky,monsky2}.

\section{Combination with descent}
To find Heegner points of large height, say 500 or more, it is usually
best first to do a descent on the elliptic curve, as this will tend to
reduce the size of the rational point by a significant factor.\footnote{
 During the mid 1990s, Cremona and Siksek worked
 out a few examples using~\hbox{2-descent}.}
Upon doing a \hbox{2-descent}, we need only $H/3$ digits of precision
if we represent the covering curve as an intersection of quadrics in~$\P^3$
and use \hbox{4-dimensional} lattice reduction,
and if we do a \hbox{4-descent} we need only $H/12$ digits.
We first explain how these lattice reduction methods work,
and then show how to use them in our application.
It might also be prudent to point out that if $E$ has nontrivial
rational isogenies, then one should work with the isogenous curve
for which the height of the generator will be the smallest.\footnote{
 Because $\Sha$ might have different size for the various isogenous curves,
 we cannot always tell beforehand which curve(s) will have a generator of
 smallest height.}

\subsection{Lattice Reduction}
Most of the theory here is due to Elkies~\cite{elkies-ants4}.
We first describe a \hbox{$p$-adic} method --- this is not immediately
relevant to us as we do not know how to approximate the Heegner
point in such a manner, but it helps to understand the idea.
Let $F(W,X,Y)=0$ be a curve in~$\P^2$. We wish to find rational points on~$F$.
Let $(1:x_s:y_s)$ be a (nonsingular) point modulo some prime~$p$,
and lift this to a solution $(1:x_0:y_0)$ modulo~$p^2$.
Then determine $d$ such that any linear combination
of $(1:x_0:y_0)$ and $(0:p:dp)$ will be a solution mod~$p^2$
(computing~$d$ essentially involves taking a derivative).
Then perform lattice reduction on the rows of the matrix
$$\begin{pmatrix} 1&x_0&y_0\\ 0&p&dp\\ 0&0&p^2 \end{pmatrix}.$$
Finally search for global solutions to~$F$ by taking small
linear combinations of the rows of the lattice-reduced matrix.
If we choose $p$ to be around $B$ for some height bound~$B$,
upon looping through all local solutions modulo~$p$ we should
find all global points whose coordinates are of size~$B$;
in general we take $p$ of size $B^{2/n}$ in projective~\hbox{$n$-space}.
This can be used, for instance, to search for points on a cubic model
of an elliptic curve.\footnote{
 This description is due to Elkies and
 is noted by Womack (\cite[Section 2.9]{Womack}).}

Over the real numbers the description is more complicated.
Here we deal with the transformation matrix of the lattice reduction.
If we wanted to do \hbox{2-dimensional} reduction, that is,
continued fractions, on a real number~$x_0$, we would
perform lattice reduction on the rows of the matrix
$$M_2=\begin{pmatrix} 1&-x_0B\\ 0&B \end{pmatrix}$$
to get good rational approximations to~$x_0$.
We can note that $\overrightarrow{(1,x_0)}M_2=\overrightarrow{(1,0)}$
and that the transformation matrix~$T$ for which $TM_2$ is lattice-reduced
has the property that $\overrightarrow{(1,0})T$ is approximately
proportional to~$\overrightarrow{(1,x_0)}$.
In four dimensions we take a point $(1:x_0:y_0:z_0)$ on some curve,
assuming that derivatives of $y$ and $z$ with respect to~$x$
are defined at this point. The matrix we use here is\footnote{
 The \hbox{3-dimensional} version is just the upper-left corner.}
$$M_4=\begin{pmatrix}
1&-x_0B&(y'x_0-y_0)B^2&\bigl(-e(y'x_0-y_0)+z'x_0-z_0\bigr)B^3\\
0&B&y'B^2&(ey'-z')B^3\\ 0&0&B^2&-eB^3\\ 0&0&0&B^3 \end{pmatrix}.$$
Here $e=z''/y''$ and all the derivatives are with respect to~$x$
and are to be evaluated at $(1:x_0:y_0:z_0)$.
Note that if we have computed $(1:x_0:y_0:z_0)$ to $H$ digits of precision,
we must ``lift'' it to precision~$3H$ to use this.
Similar to the \hbox{2-dimensional} case of above we have that
$\overrightarrow{(1,x_0,y_0,z_0)}M_4=\overrightarrow{(1,0,0,0)}$,
and $\overrightarrow{(1,0,0,0)}T$ is approximately proportional
to~$\overrightarrow{(1,x_0,y_0,z_0)}$.

\subsection{Results}
We now combine descent with the Heegner point method.
We assume that we have a cover $C\rightarrow E$, and for each
point $\WP(\dot z)$ given by the above algorithm we compute its
real pre-images on~$C$. For a \hbox{2-covering}~quartic,
the \hbox{$x$-coordinate} has size~$H/4$,
but the \hbox{$y$-coordinate} on the quartic will be of size~$H/2$.
Either continued fractions on the~\hbox{$x$-coordinate}
or \hbox{3-dimensional} lattice reduction on both coordinates and the curve
requires a precision of $H/2$ digits --- however, if
our \hbox{2-cover} is given as an intersection of quadrics in~$\P^3$,
then we only need a precision of~$(H/2)(2/3)$ since the Elkies method
does better in higher dimension. For a \hbox{4-cover} represented as
an intersection of quadrics in~$\P^3$, the coordinates are of size~$H/8$,
and so we need a precision of $(H/8)(2/3)$ to recover our point.

We now give two examples of Heegner points of large height.\footnote{
 These examples exemplify the experimental and heuristic correlation
 between large heights and large cancellation in $c_4^3-c_6^2=1728\Delta$,
 since $\Omega_{\rm re}$ can then be unusually small for a given $|\Delta|$.}
First we consider $E$ given by $[0,1,1,-4912150272,-132513750628709]$,
for which \hbox{$N=421859$}.
Here the Heegner point is of height $3239.048+$.
We refer the reader to \cite{Womack} for how to do a \hbox{4-descent}.
The intersection of quadrics that gives the \hbox{4-cover} is
given by the symmetric matrices
$$\begin{pmatrix}1&3&14&4\\ 3&7&9&8\\ 14&9&-8&19\\ 4&8&19&13\end{pmatrix}
\qquad{\rm and}\qquad
\begin{pmatrix}
16&-10&5&-5\\ -10&29&-3&-5\\ 5&-3&-1&8\\ -5&-5&8&13\end{pmatrix}.$$
We used $D=-795$ for which the index is~4. The class group has size~4;
upon using the pairing from complex conjugation we need only 2 forms,
which we can take to be $(421859,234525,32595)$ and~$(421859,384997,87839)$.
We need to use about $3239/12\log(10)\approx 120$ digits of precision
and take around 1.3 million terms of the~\hbox{$L$-series}.
For our approximation to a generator on $\C/\Lambda$ we get
\vskip0pt\noindent
$$\dot z=\hbox{\small
0.00825831518406814312450985646222558391095207954623175715662897127635126006560626891914983130574212343000780426018430276055.}$$
\vskip0pt\noindent
We find the real pre-images of this on the \hbox{4-cover}
(this can be done via a resultant computation)
and then via \hbox{4-dimensional} lattice reduction
we obtain the point
\vskip-4pt\noindent
$$\hbox{\realsmall
(90585849222350621011339302424932326542192474474854331313031216338204053880670077944701302491852572823731202634266219944146702509489824532529044887987947859472355124939471295729,}$$
\vskip-10pt\noindent
$$\hbox{\realsmall
58207848469468567249250100904745604517491584621654101065649337689496036041318068019646159331652386264579879746727095434743171075008134671399964513924870607157340785327661071757,}$$
\vskip-10pt\noindent
$$\hbox{\realsmall
-52660183473004831875084410642918250532458007522956523515279464248174730240287018424148701587123000269221530600465998289112443045957965377530412500941202059568743656911281766881,}$$
\vskip-10pt\noindent
$$\hbox{\realsmall
120566343955994724443268162651063750286159426913472902595469949397024821835521392662156637479278825638673289873881696613629128450397637909394604102581491799075589326751709650806),}$$
which can then be mapped back to~$E$.
Even though we only used 120 digits of precision in our computation
of a real approximation to the Heegner point, we can recover a point
with approximately 3/2 as many digits.
Note that if we did not use descent, but recovered the point
on the original curve using the Cremona-Silverman method,
we would need 12 times the precision and 12 times as many terms
in the \hbox{$L$-series} --- this could be a total time factor of
as much as~$12^3$, depending upon the efficiency of our high-precision
arithmetic. This computation (including 2-descent and 4-descent which each
take a second) takes less than a minute.

Finally we give a more extreme example --- this is the largest
example which we have computed. The curve is from the database of
Stein and Watkins~\cite{SteinWatkins}.
Let $E$ be given by $[0,0,1,-5115523309,-140826120488927]$,
for which $N=66157667$ and the Heegner point is of height $12557+$.
The intersection of quadrics that gives the \hbox{4-cover} is
given by the two symmetric matrices
$$\begin{pmatrix}0&1&3&3\\ 1&5&-1&-6\\ 3&-1&8&-2\\ 3&-6&-2&16\end{pmatrix}
\qquad{\rm and}\qquad
\begin{pmatrix}
12&-21&-10&-68\\ -21&-13&-7&-27\\ -10&-7&3&-7\\ -68&-27&-7&15\end{pmatrix}.$$
We select $D=-1435$ for which our 2 forms are
$(66157667,2599591,25537)$ and $(66157667,37610323,5345323)$
and the index is~2.
We need 460 digits of precision and 600 million terms of the $L$-series.
This takes less than a day. We list the \hbox{$x$-coordinate} of the
point on the original elliptic curve. It has numerator
\def\myskip{-7.5pt}
\goodbreak\noindent
\vskip0pt\noindent
\hbox{\notsosmall
3677705371866775066140056423418271700879322694922855847262187700616535463492710
1580536513437032674306114130646450005288670465199839976647884079191530786174150
}\vskip\myskip\noindent
\hbox{\notsosmall
7273933802628157325092479708268760217101755385871816780548765478502284415627682
8471927526818990949626599378706300367603592935770218062374839710749312284163465
}\vskip\myskip\noindent
\hbox{\notsosmall
0785238169688322765007203996448159721599599329974493411710628985038936400655249
7835877740257534533113775202882210048356163645919345794812074571029660897173224
}\vskip\myskip\noindent
\hbox{\notsosmall
3703377010561657350085906402970902987091215062666972664619932018253973699995508
6814229431275632217741073053282806475960497536924235099356803072693704991160726
}\vskip\myskip\noindent
\hbox{\notsosmall
4109782746847951283794119298941214490794330902986582991229569401523519938742746
3761071907702040105138183490127866378892547110594555551738109049119276198990318
}\vskip\myskip\noindent
\hbox{\notsosmall
5514929232533858983197973702640271104974259411600038060148083998297555750603585
1728035645241044229165029649347049289119188596869401159325131363345962579503132
}\vskip\myskip\noindent
\hbox{\notsosmall
3398472754224400945538247051892256536774595128631179117218385529343091245081344
9336643740809392436203974991190741697350414232211175705858420072502263211616472
}\vskip\myskip\noindent
\hbox{\notsosmall
0164998641729522677460525999499077942125820428879526063735692685991018516862938
7960475973239865371541712483169437963732171919939969937146546295368843960579247
}\vskip\myskip\noindent
\hbox{\notsosmall
9093864765666328159617814572211609821650093033382432180672693701819013619055657
3208807048355335567078793126656928657859036779350593274598717379730880724034301
}\vskip\myskip\noindent
\hbox{\notsosmall
8677394437498418094567158841937203289014615526598826284058422097567571678166621
3994508186464210853359598997571625925924015283405094065447961714768592250085694
}\vskip\myskip\noindent
\hbox{\notsosmall
4449822045386092122409096978544817218847897640513477806598329177604246380812377
7390491844755507773416209859765703930378802827649670195524084007307548226764414
}\vskip\myskip\noindent
\hbox{\notsosmall
8171538534400197983223265241488833586556737721436045600329696166817748194480906
6257442596772347829664126972931904101685281128944780074646796760942430959617022
}\vskip\myskip\noindent
\hbox{\notsosmall
2574798740894035649650388853798178669200489298145202684936775070590737659026716
3808736648849670283632626857459312324510742034887810176312389334765702027559124
}\vskip\myskip\noindent
\hbox{\notsosmall
8824247800594270862052082185973393290009189867677259458080676065098703453539525
5769756395437005076256407298723407894063143944684005844559206833619762001218344
}\vskip\myskip\noindent
\hbox{\notsosmall
3075123390147322849749056199807848625107499352887131879740334808737042690099755
6442577081254910572185107856605139877331015042842121106080690743578173268489400
}\vskip\myskip\noindent
\hbox{\notsosmall
4990568983126219539479670123584145477528970810970917957442036976840460662556632
0124229276012675987126600451637743296191727204021714708356339998761242059527579
}\vskip\myskip\noindent
\hbox{\notsosmall
2033855676991823368254862159558450043808051481533297270035287382247038279293223
9463850701180823069589872686033969240544031038574440588486055874154005176700326
}\vskip\myskip\noindent
\hbox{\notsosmall
3112120612773248134039108827779648854441573815655301476840624615466600513969042
8085145098272500791416214774673484501826722500527091164944262537169595848931680
}\vskip\myskip\noindent
\hbox{\notsosmall
7540967747128604905727462240940311870432045261072392010796034682975228951065985
6743701508334879787536416279769396881980413954888575128268715223707826035870523
}\vskip\myskip\noindent
\hbox{\notsosmall
0284426203064493684250614282879918107733796207067250003823959412935677624093236
0470386373655773263995890088045077860119731559277310730347065365574614438066227
}\vskip\myskip\noindent
\hbox{\notsosmall
0762241108780937187215721045683689249361383679202676182038221716548199892412360
4782787923229739171920575447007099501678380795077013113325989801385729993920818
}\vskip\myskip\noindent
\hbox{\notsosmall
3016544242513395646068768201219283722462133998592132827925111680439534438397939
0113997419447930029756609766453919938465190843618873242881837330238304638885942
}\vskip\myskip\noindent
\hbox{\notsosmall
7937893841888014266685177616605644783704135794931830750265686335934066565240944
0494482130055919971289855607602603992142786359126343515867623548693540215307461
}\vskip\myskip\noindent
\hbox{\notsosmall
8999289958255459763210830963856929696480004698307273623848314901471460089605655
2029642747991419063454749142059564274298254654925893866404955146903330024475746
}\vskip\myskip\noindent
\hbox{\notsosmall
1635437149962496524201711710542317263364935415869714317789440514810596337383994
1141857432381177094972972684361267292500063135565983416420055444131545100343345
}\vskip\myskip\noindent
\hbox{\notsosmall
2466204707123811663623662837296862948061758759928631763661985185615801886205770
7210320063041448677873470583163922956715800916558720872094859132869301288586404
}\vskip\myskip\noindent
\hbox{\notsosmall
4258912545426858039748457192101231887231162489831761560762817646009744133632354
9031828235965636277950827328087547939511112374216436584203379248450122647406094
}\vskip\myskip\noindent
\hbox{\notsosmall
0351711307406637235476759398859593638811358930351020183894442127461462503283482
4261067352402237899497839202009881472197450206269281573668922975906582209394279
}\vskip\myskip\noindent
\hbox{\notsosmall
5318705345275598989426335235935505605311411301560321192269430861733743544402908
5864973053536009094312149332025225287171092144929593300160658102876231441792884
}\vskip\myskip\noindent
\hbox{\notsosmall
6666488854062270234670421375245637257444956397921578240656693788535294587199454
1770838871930542220307771671498466518108722622109421676741544945695403509866953
}\vskip\myskip\noindent
\hbox{\notsosmall
1672776282802324648392150034740488969680375446600297557400655812701390832499032
1257223041794224979546710070039394431032500967717918210997094334680733501444683
}\vskip\myskip\noindent
\hbox{\notsosmall
9612282508824324073679584122851208360459166315484891952299449340025896509298935
9393577217235439331087432419973874470183959253201676376403284079570698454395013
}\vskip\myskip\noindent
\hbox{\notsosmall
8123460586749500340201672462640085536963652115500914717624590414906922543864692
8549072337653348704931901764847439772432025275648964681387210234070849306330191
}\vskip\myskip\noindent
\hbox{\notsosmall
7903804123961154462408325834813663721323008490608352621368323153110529033675038
5743792050893130528314337942393060136915457253067727886206663888425022179164712
}\vskip\myskip\noindent
\hbox{\notsosmall
3563828956462530983567929499493346622977494903591722345188975062941907415400740
881}
\goodbreak\vskip0pt\noindent
and denominator (which is square)
\vskip-2pt\noindent
\hbox{\notsosmall
4255044272974739888318147243888463149102796676262869844129037001939069082536978
3206694892509489159845392479060050865239326960543711865223880415469150998005445
}\vskip\myskip\noindent
\hbox{\notsosmall
0063276671301195356576897799687302838382126577552634434635712354111258918875514
7040471997532333827695072210398747506143870850088480631610375575105629402969923
}\vskip\myskip\noindent
\hbox{\notsosmall
3799066968836324267266407778292384246876573476455250148990478999518159066382588
2826884437614702086601072912460205046718987111919759741772924082107906614768188
}\vskip\myskip\noindent
\hbox{\notsosmall
2592877855973838960678143778693770394049052970155019875224052204788541398018355
6155907460046868073063343753429797282371426595980125377117083031923094847501187
}\vskip\myskip\noindent
\hbox{\notsosmall
8153976528034885787379890174882324346380444041171547892001495923651535208922165
7498291013202443206177483221913021678256235579893414436481218449503035888524081
}\vskip\myskip\noindent
\hbox{\notsosmall
1574061954766438356661270643806422096593442057480569818796530004016888734707062
4862281490478227347630906482652560940554823360522517964582349779334663133777009
}\vskip\myskip\noindent
\hbox{\notsosmall
3238960301538846226719138651429924452362996719692807033817042217290429195445858
8735665717672690912582485100989972418433868929605913487922691271558721181323613
}\vskip\myskip\noindent
\hbox{\notsosmall
7001205209064857537600160042225772832307591236986881863116245504382548723553194
1123620115037512323584244014680824673299882453235229917216799756372612183353210
}\vskip\myskip\noindent
\hbox{\notsosmall
2368619557885419527341455580068205352191453004482040011648102778922881052572304
7625080568425427769667815771926774952742389468234144023574739674218399717724794
}\vskip\myskip\noindent
\hbox{\notsosmall
6449079490990330823304497772932137791366115383255184513570505904177076177300233
5204714530513728518042273513609336859761902134080074259874956489916483283558415
}\vskip\myskip\noindent
\hbox{\notsosmall
1879023308997164082639142186544641615917653330537431550152668753631028786525822
5341831409297024404559487986223470061428509852910119956281813557240650373612303
}\vskip\myskip\noindent
\hbox{\notsosmall
8061242168264761915827597120621691107510211062002225029210397375828078654532129
3000346117555930108145915118154377991059233962745297863017035414522428825832481
}\vskip\myskip\noindent
\hbox{\notsosmall
1018569690301921451517013795632224485059295610701352935567534896297148895118420
8519602258834025155300581173742601529286315229193789791848640470464974617150056
}\vskip\myskip\noindent
\hbox{\notsosmall
6225974820050798241870400003418958513067964846039950289279541943230841117696177
7683169393789440445404586765622160449971094877274592261820143359041089932670055
}\vskip\myskip\noindent
\hbox{\notsosmall
5129814741781368814174760031579179374183703665646076025908428641450401691692865
3373373749351673138873697355678073568445462500597124999569428768348290694470131
}\vskip\myskip\noindent
\hbox{\notsosmall
3969575867184534505315069252495619031118665012176098178260784016899909892825429
4766906595754090551526543788821211073870016226804688310128261364457467501381181
}\vskip\myskip\noindent
\hbox{\notsosmall
9989545229962754945276194487111239718227554281251751692833307185347349647947117
7474588594469094399278566301150327219929702835948186371102871019478295556041049
}\vskip\myskip\noindent
\hbox{\notsosmall
9456979053009021516126703280144207697842588434057829068494470566384220387343494
6895313693410268205549042641332729214043530781277783400932521461065846050516266
}\vskip\myskip\noindent
\hbox{\notsosmall
0636930982933352215314240689644147222987920125050724733503025375041491926808614
2811991181135615102063450273874178955396056033628717223631646597545299321971795
}\vskip\myskip\noindent
\hbox{\notsosmall
7208380379846976284134301874897304106795470358899054823624849472365173904528952
1882710298950201026345472967707812057733059495444795253020394607206147343184720
}\vskip\myskip\noindent
\hbox{\notsosmall
1289138202794711371684597932965767939112179634242219380882552731959067823646242
1470515380610223708802325911166125101655807805544810736511736112015725788207983
}\vskip\myskip\noindent
\hbox{\notsosmall
0434627788427533485656242565476771106266450533897646444735950588164991426098055
5651233453080870411547161166201727321046199497573313805291689786242134493722493
}\vskip\myskip\noindent
\hbox{\notsosmall
6540132813720819087960384794865152730525579241209808589941796218972732173245666
6040832323196614881188508752637063922635721727888588861006697067295082662856621
}\vskip\myskip\noindent
\hbox{\notsosmall
4278955509812709830759473382274095888608048326897787309857569462281336128176938
4832759774254005221673936338256702877699726892622342424638164907882484814569895
}\vskip\myskip\noindent
\hbox{\notsosmall
7882084132643790230707482558729729822722846435295235634841524505128041497093453
4815003363838566586486090914244124413024752619262170944161605020754084937204942
}\vskip\myskip\noindent
\hbox{\notsosmall
1637232671248907772397162458935318730166463274525319193925966950757860571076838
7023661969449344902879751178806889386035847044662250758843856626244205371560027
}\vskip\myskip\noindent
\hbox{\notsosmall
8734508668324795740959080113149044555394947297908982119899939239220581201178632
4196553409374966862790769096287105922007091926687147181741236537543028535996580
}\vskip\myskip\noindent
\hbox{\notsosmall
4594745527651369886073387194293961511036178306538077683719722004972020545664652
0373809780940384932684877526623946791294274354580233914301360048309031821084423
}\vskip\myskip\noindent
\hbox{\notsosmall
7242526421831575726183780763818586110019840381512817836756338699203654792269435
8498639362771798041618046013662235878191246173024597377861825966407677879013966
}\vskip\myskip\noindent
\hbox{\notsosmall
4482603550041332494718324185233462156771301830971108909989639075891161512882964
6468875441800478905861321031022715201049054765075321382013759982613256651365903
}\vskip\myskip\noindent
\hbox{\notsosmall
0867731020594559477466913485854520908208415598659344363758912591077035997911074
2764732875552438147802286623083638506353523937829379144982925749510228726599819
}\vskip\myskip\noindent
\hbox{\notsosmall
6492322685048518708341003752717812628525671219633053830557941194740861837064808
1044119314858092966789985880558618838749736465935723939056195292855739824197002
}\vskip\myskip\noindent
\hbox{\notsosmall
8574716638792985681864321864807876516532795271727594036621733392160501881017349
0426678520993395921957402807842113084883601771499771585515616148588419257652541
}\vskip\myskip\noindent
\hbox{\notsosmall
1096570236580006784778240115504215083155473114336654532035309712573054076300648
0842903826240579636224118544527127667643337427417482571764186299386690494751101
}\vskip\myskip\noindent
\hbox{\notsosmall
45831487629518887622732479751575475237265720522099995921350373526902203987716.
}

\end{document}